\documentclass[11pt,a4paper,fleqn]{amsart}
\usepackage{a4wide,amsfonts,amsmath,latexsym,amssymb,euscript,graphicx,units,mathrsfs}

\usepackage{graphicx}
\usepackage{color}
\usepackage{amssymb}
\usepackage{amssymb}
\usepackage[T1]{fontenc}
\usepackage{latexsym}
\usepackage{xypic}
\usepackage{eufrak}
\usepackage{euscript}
\usepackage{amsfonts,amsmath}
\usepackage{verbatim}
\usepackage{fancyhdr}
\usepackage[english]{babel}
\usepackage{mathrsfs}
\usepackage{units}

\newtheorem{prop}{Proposition}[section]
\newtheorem{theorem*}{Theorem}
\newtheorem{cor}[prop]{Corollary}
\newtheorem{lemma}[prop]{Lemma}
\newtheorem{lemme}[prop]{Lemma}

\newtheorem{theorem}[prop]{Theorem}

\newtheorem{defi}[prop]{Definition}

\newtheorem*{prop-cont}{Proposition \ref{azerty}}
\newtheorem*{theorem-cont}{Theorem \ref{meckes}}

\renewcommand{\geq}{\geqslant}
\def\leq{\leqslant}

\newcommand{\R}{\mathbb{R}}

\def\1{{\mathbf{1}}}

\def\1{{\mathbf{1}}}
\def\0.5{{\frac{1}{2}}}

\def\E{\mathbb{E}}

\begin{document}

\begin{center}
{\Large{\bf Stein's method for asymmetric $\alpha$-stable distributions, \\
with application to the stable CLT}}\\~\\
Peng Chen, Ivan Nourdin and Lihu Xu \\
{\it University of Macau}, {\it Universit\'e du Luxembourg} and {\it University of Macau}\\~\\
\end{center}
{\small \noindent {\bf Abstract:}
This paper is concerned with the Stein's method associated with a
(possibly) asymmetric $\alpha$-stable distribution $Z$, in dimension one.
More precisely, its goal is twofold.
In the first part, we exhibit a genuine bound for the Wasserstein distance between $Z$ and any integrable random variable $X$, in terms of an operator that reduces to the classical fractional Laplacian  in the symmetric case.
Then, in the second part we apply the aforementioned bound to compute
error rates in the stable central limit theorem, when the entries are
in the domain $\mathcal{D}_\alpha$ of normal attraction of a stable law of exponent $\alpha$. To conclude, we study the specific case where the entries are Pareto like multiplied by a slowly varying function, which provides an example of random variables that do not belong to $\mathcal{D}_\alpha$, but for which
our approach continues to apply.
 \\


\bigskip

\section{Motivation and main results}

Sums of independent and identically distributed random variables are among the most basic quantities we can encounter in probability theory. When the entries are square integrable, in the large limit they become approximately gaussian, thanks to the usual central limit theorem. And provided we assume a little more integrability for the entries, we can also compute explicit error rates by means of {\it Fourier analysis} (see, e.g., Berry-Esseen theorem and relatives).

What is nowadays referred to as {\it Stein's method} is another technique, invented by the great statistician Charles Stein in the late sixties, to compute explicit rates for the error in the gaussian approximation. The power of Stein's method compared to Fourier analysis is that one can also use the former in the presence of dependence (see \cite{CGS} for a self-contained treatment of Stein's method and its many ramifications).
Over the years, Stein's method has become an indispensable tool in probability theory and statistics, with applications in a bunch of different areas.
Initially proposed for gaussian approximation, there now exist
extensions in many other contexts, ranging from classical stochastic approximations such as poisson or gamma, to more exotic ones.
For a regularly updated list of all available extensions in the literature, we refer the reader to the very useful webpage \cite{yvik-web} maintained by Yvik Swan.

Surprisingly, and despite the fact that  the \emph{stable central limit theorem} is undoubtedly among the most important limit theorems in probability theory,
Stein's method for stable approximation has been barely developed.
Actually, to the best of our knowledge only one paper  (\cite{lihu}) is concerned with this problem, whereas two others \cite{AH,AMPS} develop the so-called Stein-Tikhomirov approach (a kind of mixture of Stein's method and Fourier analysis).
The main difference between \cite{lihu} and the present paper
is that the former only considers the symmetric case and develops the $K$-function method to compute bounds, whereas here we consider more generally the asymmetric case and we shall develop a variant of the leave-one out approach to compute our bounds.

To describe our results in a more explicit way, we start by recalling
the definition of an $\alpha$-stable distribution.
Note that we only consider the case $\alpha\in(1,2)$ starting from now.
This is because it is the only range
for $\alpha$ that may make sense when we work with the Wasserstein distance
$d_W$ (defined, for two integrable random variables $X$ and $Y$, as
$d_W(\mathcal{L}(X),\mathcal{L}(Y))=\sup \big| \E[h(X)]-\E[h(Y)]\big|$,
where the supremum runs over all $1$-Lipschitz functions $h:\R\to\R$). In addition, we denote $C^{k}(\R)=\{f:$ $\mathbb{R}\rightarrow\mathbb{R};$ $f,~f',~\cdots,~f^{(k)}$ are all continuous functions$\}.$
\begin{defi}\label{defstable}
Let $\alpha\in(1,2)$, $\sigma> 0$ and $\beta\in [-1,1]$ be real numbers.
\begin{enumerate}
\item[1a)] We say that $Z$ is distributed according to the \emph{$\alpha$-stable law with parameters $\sigma$ and $\beta$}, and we write $Z\sim S_\alpha(\sigma,\beta)$, to indicate that
\[
\E[e^{itZ}]=
{\rm exp}\big\{-\sigma^\alpha|t|^\alpha(1-i\,\beta\,{\rm sign}(t)\,\tan\frac{\pi\alpha}2)\big\},\quad t\in\R.
\]
\item[1b)] In the particular case where $\beta=0$, we have
\[
\E[e^{itZ}]=
{\rm exp}\big\{-\sigma^\alpha|t|^\alpha\big\},\quad t\in\R,
\]
and we
say that $Z$ is distributed according to the \emph{symmetric $\alpha$-stable law of parameter} $\sigma$, and we write $Z\sim S\alpha S(\sigma)$.
\item[2)] For $\phi:\R\to\R$ in $C^2(\R)$  such that
$\|\phi''\|_\infty<\infty$, we set
\[
(\mathcal{L}^{\alpha,\beta}\phi)(x)=d_\alpha \,
\int_\R \frac{\phi(u+x)-\phi(x)-u\phi'(x)}{2|u|^{1+\alpha}}
\left[
(1+\beta){\bf 1}_{(0,\infty)}(u)+(1-\beta){\bf 1}_{(-\infty,0)}(u)
\right]
 du,
\]
where $d_\alpha = \left(
\int_0^\infty \frac{1-\cos y}{y^{1+\alpha}}dy
\right)^{-1}$.
\end{enumerate}
\end{defi}
It is immediate to check that $Z/\sigma\sim  S_\alpha(1,\beta)$  iff $Z\sim  S_\alpha(\sigma,\beta)$. This is why, starting from now and without loss of generality, we will only consider stable distributions for which $\sigma=1$.
Moreover, when $\beta=0$ we observe that $\mathcal{L}^{\alpha,0}$ reduces to the  fractional Laplacian $\Delta^\frac\alpha2$ of order $\alpha/2$, see e.g. \cite[Section 4.1]{lihu}.

It is well known that the distribution $S_{\alpha}(\sigma,\beta)$ admits a smooth density \cite[Proposition 2.5, (xii)]{sato}, denote it by  $p_{\sigma,\beta}$, which satisfies
\begin{equation}\label{def-den}
\int_\R e^{i\lambda x}p_{\sigma,\beta}(x)dx = {\rm exp}\big(-\sigma^\alpha|\lambda|^\alpha(1-i\,\beta\,{\rm sign}(\lambda)\,\tan\frac{\pi\alpha}2)\big),\quad \lambda\in\R.
\end{equation}

Our main first result is the following theorem, that provides a bound
 for the Wasserstein distance $d_W$ between $Z\sim S_\alpha(1,\beta)$ and any integrable random variable $X$, in terms of the operator $\mathcal{L}^{\alpha,\beta}$.

\begin{theorem}\label{thm1}
Let $\alpha\in(1,2)$ and $\beta\in [-1,1]$ be real numbers, and let
$Z\sim S_\alpha(1,\beta)\sim p_Z(z)dz$.
Then, for any integrable random variable $X$,
\begin{equation}\label{steinstable}
d_W(\mathcal{L}(X),\mathcal{L}(Z))
\leq \sup_{\substack{\|\phi'\|_\infty\leq \alpha\\\|\phi''\|_\infty\leq \eta_{\alpha,\beta}}} \big|\E[(\mathcal{L}^{\alpha,\beta}\phi)(X)]-\frac1\alpha \E[X\phi'(X)]\big|,
\end{equation}
where $\eta_{\alpha,\beta}={\rm Beta}(\frac2\alpha,1-\frac1\alpha)\int_\R\big|p'_Z(y)|dy<\infty$.
\end{theorem}

To illustrate a possible and explicit use of our abstract Theorem \ref{thm1},
we will now describe our second main result, which
computes rates in the stable central limit theorem by means of
a leave-one-out approach.

First, let us mention that the problem of calculating rates
in the stable central limit theorem  is of course not new, and there is a dedicated literature
on the subject dating back to the seventies (see, e.g., \cite{lithuanian} and references therein). At that time, the
main challenge was to compute
rates in the \emph{Kolmogorov distance} (written $d_{\rm Kol}$ in the sequel)
by means of Fourier analysis. A representative result obtained in this family of papers
 is the following estimate taken from \cite{banis}.
 Assume that $X_1,X_2,\ldots$ are independent copies drawn from the \emph{Pareto law of index $\alpha\in(1,2)$}, that is, suppose that the
common density is $p(x)=\frac{\alpha}{2}|x|^{-(1+\alpha)}{\bf 1}_{[1,\infty)}(|x|)$.
Then, with $\sigma=\left(\frac{\alpha}2\int_\R \frac{1-\cos y}{|y|^{1+\alpha}}dy\right)^{\frac1\alpha}$ and $S_n=\frac{1}\sigma n^{-\frac1\alpha}(X_1+\ldots+X_n)$, one has
\begin{equation}
d_{\rm Kol}(\mathcal{L}(S_n),S\alpha S(1))
=O(n^{1-\frac2\alpha})\quad\mbox{as $n\to\infty$.}\label{dKol}
\end{equation}

This being recalled, let us now give the estimates we obtain on our side from Theorem \ref{thm1}, after having combined it with a variant of the \emph{leave-one-out approach}
classically used in the context of normal approximation (see Section \ref{loo}).
Since we rely on Theorem \ref{thm1}, our results are given in terms of the
Wasserstein distance $d_W$ rather than the Kolmogorov distance $d_{\rm Kol}$, as in (\ref{dKol}); 
note that there is no subordination relationship between $d_{\rm Kol}$ and $d_W$, so  bounds in either metrics are relevant in their own right.

Before giving the second main theorem, we first recall the definition of normal attraction of a stable law of exponent $\alpha$.

\begin{defi}
If $X$ has a distribution function of the form
\begin{equation}\label{px}
F_X(x)=\big(1-\frac{A+\epsilon(x)}{|x|^{\alpha}}(1+\beta)\big){\bf 1}_{[0,\infty)}(x)+\frac{A+\epsilon(x)}{|x|^{\alpha}}(1-\beta){\bf 1}_{(-\infty,0)}(x)\big),
\end{equation}
where $\alpha\in(1,2)$, $A>0$, $\beta\in[-1,1]$ and  $\epsilon: \R\to\R$ is a bounded function vanishing at $\pm\infty$, then we say that $X$ is in the domain $\mathcal{D}_\alpha$ of normal attraction of a stable law of exponent $\alpha$.
\end{defi}

In (\ref{px}), the function $\epsilon$ is supposed to be bounded, that is, there exists $K>0$ such that $|\epsilon(x)|\leq K$. More specifically,
let us assume the existence of $K>0$ and $\gamma\geq 0$ such that
\begin{align}\label{bound}
|\epsilon(x)|\leq\frac{K}{|x|^{\gamma}},\quad x\neq0.
\end{align}

Observe that making $\gamma=0$ in (\ref{bound}) simply means that we do not want to make any extra assumption on $\epsilon$ defined in (\ref{px}).

We can now state our second main theorem.
\begin{theorem}\label{thm2}
Let $X_1,X_2,\ldots$ be independent and identically distributed random variables defined on a common probability space, and suppose that $X_1$ has
a distribution of the form (\ref{px}) with $\epsilon(x)$ satisfying (\ref{bound}).
Set
\begin{equation}\label{sn}
S_n=\frac{1}{\sigma}n^{-\frac{1}{\alpha}}(X_1+\ldots+X_n-n\E[X_1]),
\end{equation}
where
$\sigma=\left(A\alpha\int_\R \frac{1-\cos y}{|y|^{1+\alpha}}dy\right)^\frac1\alpha$.
Then there exists $c_{\alpha,\beta,\gamma}$ (that can be made explicit) depending only on $\alpha$, $\beta$ and $\gamma$ such that,
\begin{align*}
d_W(\mathcal{L}(S_n),S_\alpha (1,\beta))
\leq&c_{\alpha,\beta,\gamma}
\begin{cases}
n^{1-\frac{2}{\alpha}},\quad &\gamma\in(2-\alpha,\infty),\\
n^{1-\frac{2}{\alpha}}\log n,\quad &\gamma=2-\alpha,\\
n^{-\frac{(\alpha-1)\gamma}{\alpha(1-\gamma)}},\quad &\gamma\in(0,2-\alpha),\\
n^{1-\frac{2}{\alpha}}\int_{-\sigma n^{\frac{1}{\alpha}}}^{\sigma n^{\frac{1}{\alpha}}}\frac{|\epsilon(x)|}{|x|^{\alpha-1}}dx+(\sup_{|x|\geq \sigma n^{\frac{1}{\alpha}}}|\epsilon(x)|)^{\alpha-1}, &\gamma=0.
\end{cases}
\end{align*}
\end{theorem}

Since, by Ces\'aro, the quantity $n^{1-\frac2\alpha}\int_{-\sigma\,n^{\frac1\alpha}}^{\sigma\,n^{\frac1\alpha}}
\frac{|\epsilon(x)|}{|x|^{\alpha-1}}dx$  tends to zero when $\epsilon(x)\to 0$ as $x\to\pm\infty$,
it is immediate to deduce from Theorem \ref{thm2} the following corollary,
which is of independent interest.
\begin{cor}
Keep the same notation and assumptions as in Theorem \ref{thm2}.
Then, as $n\to\infty$,
$d_W(\mathcal{L}(S_n),S_\alpha (1,\beta))\to 0$.
\end{cor}

Our next result gives an improved upper bound on $d_W(\mathcal{L}(S_n),S_\alpha (1,\beta))$ under slightly more restrictive conditions (see, e.g., \cite[Theorem 2]{hall}).
\begin{theorem}\label{thm3}
Keep the same notation and assumptions as in Theorem \ref{thm2}. In addition, we further assume $\frac{\epsilon(x)}{x^{\alpha}}{\bf 1}_{(0,\infty)}(x)$ and $\frac{\epsilon(x)}{|x|^{\alpha}}{\bf 1}_{(-\infty,0)}(x)$ are ultimately monotone (that is, there exist $x_{0}>0$ such that $\frac{\epsilon(x)}{x^{\alpha}}{\bf 1}_{(0,\infty)}(x)$ and $\frac{\epsilon(x)}{|x|^{\alpha}}{\bf 1}_{(-\infty,0)}(x)$ are monotone for any $|x|>x_{0}$). Then there exists $\hat{c}_{\alpha,\beta}$ (that can be made explicit) depending only on $\alpha$ and $\beta$ such that,
\begin{equation*}
d_W(\mathcal{L}(S_n),S_\alpha (1,\beta))
\leq \hat{c}_{\alpha,\beta}\left(
n^{1-\frac2\alpha}\left[1+
\int_{-\sigma\,n^{\frac1\alpha}}^{\sigma\,n^{\frac1\alpha}}
\frac{|\epsilon(x)|}{|x|^{\alpha-1}}dx\right]
+\sup_{|x|\geq \sigma\,n^{\frac1\alpha}}|\epsilon(x)|
\right).
\end{equation*}
\end{theorem}

\medskip

To conclude this introduction, let us analyse the conclusion of Theorem \ref{thm2} and Theorem \ref{thm3} in several specific and explicit examples.

\medskip

{\it \underline{Example 1}: Pareto}.
Our first example is the simplest possible situation, that is, the case where $X_1$ is distributed according to a (possibly non-symmetric)
Pareto distribution of the form
\begin{align*}
\mathbb{P}(X_{1}> x)=\frac{1+\beta}{2|x|^{\alpha}},\quad x\geq1,\qquad \mathbb{P}(X_{1}\leq x)=\frac{1-\beta}{2|x|^{\alpha}},\quad x\leq-1,
\end{align*}
with $\alpha\in(1,2)$ and $\beta\in[-1,1]$;
in this case, (\ref{px}) holds with $A=\frac{1}2$ and $\epsilon(x)=-\frac{1}2{\bf 1}_{(-1,1)}(x)$, and we deduce from Theorem \ref{thm2} with $\gamma=0, K=\frac{1}{2}$ that
$d_W(\mathcal{L}(S_n),S_\alpha(1,\beta))=O(n^{1-\frac2\alpha})$,
compare with (\ref{dKol}).

\medskip

{\it \underline{Example 2}: Sum of two power functions of different orders}. We consider
\begin{align*}
&\mathbb{P}(X_{1}> x)=(A|x|^{-\alpha}+\widetilde{A}|x|^{-\widetilde{\alpha}})(1+\beta),\quad x\geq1,\\
&\mathbb{P}(X_{1}\leq x)=(A|x|^{-\alpha}+\widetilde{A}|x|^{-\widetilde{\alpha}})(1-\beta),\quad x\leq-1,
\end{align*}
with $1<\alpha<2$, $\alpha<\widetilde{\alpha}$, $A+\widetilde{A}=\frac12$ and $\beta\in[-1,1]$;
in this case, (\ref{px}) holds with
\[
\epsilon(x)=\widetilde{A}|x|^{\alpha-\widetilde{\alpha}}{\bf 1}_{[1,\infty)}(|x|)-A{\bf 1}_{(-1,1)}(x),
\]
and we deduce from Theorem \ref{thm3} that
$d_W(\mathcal{L}(S_n),S_\alpha(1,\beta))=\left\{
\begin{array}{lll}
O(n^{1-\frac{\widetilde{\alpha}}\alpha})&\mbox{if $\widetilde{\alpha}<2$}\\
O(n^{1-\frac{2}\alpha}\log n)&\mbox{if $\widetilde{\alpha}=2$}\\
O(n^{1-\frac{2}\alpha})&\mbox{if $\widetilde{\alpha}>2$}
\end{array}
\right..
$

\medskip

{\it \underline{Example 3}: Sum of two power functions of different orders with trigonometric function}. We consider
\begin{align*}
&\mathbb{P}(X_{1}> x)=|x|^{-\alpha}\left(A+\frac{B\sin x}{|x|}\right)
(1+\beta),\quad x\geq 1,\\
&\mathbb{P}(X_{1}\leq x)=|x|^{-\alpha}\left(A+\frac{B\sin x}{|x|}\right)
(1-\beta),\quad x\leq-1,
\end{align*}
with $1<\alpha<2$, $A+B\sin1=\frac12$ and $\beta\in[-1,1]$;
in this case, (\ref{px}) holds with
\[
\epsilon(x)=\frac{B\sin x}{|x|}{\bf 1}_{[1,\infty)}(|x|)-A{\bf 1}_{(-1,1)}(x),
\]
and we deduce from Theorem \ref{thm2} with $\gamma=1, K=B$ that
$d_W(\mathcal{L}(S_n),S_\alpha(1,\beta))=O(n^{1-\frac{2}\alpha}).$
\medskip

{\it \underline{Example 4}: Pareto with modified logarithmic tails}. We consider
\begin{align*}
&\mathbb{P}(X_{1}> x)=|x|^{-\alpha}\left(A+\frac{D}{\log|x|}\right)
(1+\beta),\quad x\geq e,\\
&\mathbb{P}(X_{1}\leq x)=|x|^{-\alpha}\left(A+\frac{D}{\log|x|}\right)
(1-\beta),\quad x\leq-e,
\end{align*}
with $\alpha\in(1,2)$, $\beta\in[-1,1]$ and
suitable $A$ and $D$;
in this case, (\ref{px}) holds with
\[
\epsilon(x)=\frac{B}{\log |x|} {\bf 1}_{[e,\infty)}(|x|)-A{\bf 1}_{(-e,e)}(x),
\]
and we deduce from Theorem \ref{thm3} that
$d_W(\mathcal{L}(S_n),S_\alpha(1,\beta))=O(\frac1{\log n})$.

\medskip
Each of the previous four examples leads to a function $\epsilon$ that satisfies
$\epsilon(x)\to 0$ as $x\to\pm\infty$, as is required in Theorem \ref{thm2} and Theorem \ref{thm3}.
But that $\epsilon$ vanishes is not a necessary condition for the stable CLT to hold.
Actually, by slightly modifying the approach leading to Theorem \ref{thm2} we can also consider examples where $\epsilon$ is a slowly varying function diverging at infinity.
Because it would be too technical to state such result
at a great level of generality,
we prefer to illustrate an explicit situation for
which our methodology still allows to conclude. Here we give a simpler proof that rather relies on the density function; note however that it would have been equally difficult to deal with the distribution function instead.
\medskip

{\it \underline{Example 5}: Pareto multiplied by a slowly varying function}. We consider
\[
p_X(x)=\frac{\alpha^2e^\alpha}{2(1+\alpha)}\,\frac{\log |x|}{|x|^{\alpha+1}}{\bf 1}_{[e,\infty)}(|x|), \quad\mbox{with $\alpha\in(1,2)$}.
\]
For the partial sums  $S_n$ to converge to the symmetric $\alpha$-stable distribution,
we must modify the normalization given in (\ref{sn}) (observe that
$\E[X_1]=0$ here). Define
the sequence $(\gamma_n)_{n\geq 1}$ implicitly by $\gamma_n=\big(n\log \gamma_n\big)^\frac1\alpha$
and set $\sigma =\left(\frac{\alpha^2e^\alpha}{(1+\alpha)d_\alpha}\right)^{\frac1\alpha}$.
We can deduce from a suitable modification of Theorem \ref{thm2} (see Section 4)  that
\begin{equation}\label{ex4}
d_W(\mathcal{L}(\frac1{\sigma\,\gamma_n}(X_1+\ldots+X_n)),S\alpha S(1))=O\big((\log n)^{-1}\big).
\end{equation}

\medskip

The rest of the paper is organized as follows.
In Section 2 we give useful properties of the operator $\mathcal{L}^{\alpha,\beta}$, we study the Stein's equation for asymmetric $\alpha$-stable distributions, and eventually we make the proof of Theorem \ref{thm1}.
In Section 3, we develop the leave-one-out approach associated with the $\alpha$-stable distribution; then, we provide the proof of Theorem \ref{thm2} and Theorem \ref{thm3}.
Finally, Section 4 is devoted to the proof of (\ref{ex4}).
\section{Proof of Theorem \ref{thm1}}

\subsection{About the operator $\mathcal{L}^{\alpha,\beta}$}

The following proposition gathers  useful alternate expressions for the operator $\mathcal{L}^{\alpha,\beta}$ introduced in Definition \ref{defstable}.

\begin{prop}\label{properties}
Fix $\alpha\in(1,2)$ and $\beta\in [-1,1]$.
Let $\phi\in C^2(\R)$ be such that $\|\phi''\|_\infty<\infty$.
We have, for all $x\in\R$ and $a>0$,
\begin{eqnarray}
(\mathcal{L}^{\alpha,\beta}\phi)(x)&=&\frac{d_\alpha}{\alpha}\int_0^\infty
\frac{(1+\beta)\big(\phi'(x+u)-\phi'(x)\big)-(1-\beta)\big(\phi'(x-u)-\phi'(x)\big)}{2u^\alpha}du\notag\\
&=& a^{1-\alpha}\frac{d_\alpha}{\alpha}\int_\R u\big(\phi'(x+au)-\phi'(x)\big)
\frac{(1+\beta){\bf 1}_{(0,\infty)}(u)+(1-\beta){\bf 1}_{(-\infty,0)}(u)}
{2|u|^{\alpha+1}}du.\label{crucial}
\end{eqnarray}
\end{prop}
\noindent
{\it Proof}. 1. One can write
\begin{eqnarray*}
&&\frac1{d_\alpha}(\mathcal{L}^{\alpha,\beta}\phi)(x)\\
&=&(1+\beta)\int_0^\infty \frac{du}{2u^{1+\alpha}}\int_0^udt\,\big(\phi'(x+t)-\phi'(x)\big)
-(1-\beta)\int_{-\infty}^0  \frac{du}{(-u)^{1+\alpha}}\int_u^0 dt\,\big(\phi'(x+t)-\phi'(x)\big)\\
&=&(1+\beta)
\int_0^\infty dt\,\big(\phi'(x+t)-\phi'(x)\big)\int_t^\infty \frac{du}{2u^{1+\alpha}}
-(1-\beta)
\int_{-\infty}^0 dt\,\big(\phi'(x+t)-\phi'(x)\big)\int_{-\infty}^t \frac{du}{2(-u)^{1+\alpha}}\\
&=&\frac{1}{\alpha}\int_0^\infty
\frac{(1+\beta)\big(\phi'(x+t)-\phi'(x)\big)-(1-\beta)\big(\phi'(x-t)-\phi'(x)\big)}{2t^\alpha}dt.
\end{eqnarray*}
2. One can write
\begin{eqnarray*}
\frac1{d_\alpha}(\mathcal{L}^{\alpha,\beta}\phi)(x)
&=&\frac{1+\beta}{\alpha}
\int_0^\infty \big(\phi'(x+t)-\phi'(x))\frac{dt}{2t^\alpha}
-\frac{1-\beta}{\alpha}
\int_{-\infty}^0 \big(\phi'(x+t)-\phi'(x))\frac{dt}{2|t|^\alpha}\\
&=&\frac{1}{\alpha}\int_\R t\big(\phi'(x+t)-\phi'(x)\big)
\frac{(1+\beta){\bf 1}_{(0,\infty)}(t)+(1-\beta){\bf 1}_{(-\infty,0)}(t)}
{2|t|^{\alpha+1}}dt\\
&=&\frac{a^{1-\alpha}}{\alpha}\int_\R u\big(\phi'(x+au)-\phi'(x)\big)
\frac{(1+\beta){\bf 1}_{(0,\infty)}(u)+(1-\beta){\bf 1}_{(-\infty,0)}(u)}
{2|u|^{\alpha+1}}du.
\qed
\end{eqnarray*}

\bigskip

Another important property of the operator $\mathcal{L}^{\alpha,\beta}$ is that it transforms
$C^2_b(\R)$-functions into $(2-\alpha)$-H\"older continuous functions, here $C^{2}_{b}(\R)=\{f\in C^{2}(\R):$ $\|f''\|_{\infty}<\infty\}.$

\begin{prop}
Fix $\alpha\in(1,2)$ and $\beta\in [-1,1]$.
Let $\phi\in C^2(\R)$ be such that $\|\phi''\|_\infty<\infty$.
We then have, for any $x,y\in\R$,
\begin{equation}
\left|
(\mathcal{L}^{\alpha,\beta}\phi)(x) - (\mathcal{L}^{\alpha,\beta}\phi)(y)
\right|\leq
\frac{4d_{\alpha}\|\phi''\|_\infty}{\alpha(2-\alpha)(\alpha-1)}\,|x-y|^{2-\alpha}.\label{crucialbis}
\end{equation}
\end{prop}
\noindent
{\it Proof}. Using (\ref{crucial}) with $a=1$ we can write, for any $x,y\in\R$:
\begin{eqnarray}
&&\left|
(\mathcal{L}^{\alpha,\beta}\phi)(x) - (\mathcal{L}^{\alpha,\beta}\phi)(y)
\right|\notag\\
&=&\left|\frac{d_\alpha}{\alpha}
\int_\R t\big(
\phi'(x+t)-\phi'(x)-\phi'(y+t)+\phi'(y)
\big) \frac{(1+\beta){\bf 1}_{(0,\infty)}(t)+(1-\beta){\bf 1}_{(-\infty,0)}(t)}{2|t|^{1+\alpha}}dt\right|\notag\\
&\leq&\frac{d_\alpha}{\alpha}
\int_\R\big|
\phi'(x+t)-\phi'(x)-\phi'(y+t)+\phi'(y)
\big| \frac{dt}{|t|^{\alpha}}\notag\\
&\leq&\frac{4d_{\alpha}}{\alpha}\|\phi''\|_\infty |x-y|\int_{|x-y|}^\infty t^{-\alpha}dt
+\frac{4d_{\alpha}}{\alpha}\|\phi''\|_\infty \int_{0}^{|x-y|} t^{1-\alpha}dt=\frac{4d_{\alpha}\|\phi''\|_\infty}{\alpha(2-\alpha)(\alpha-1)}\,|x-y|^{2-\alpha}.\notag
\end{eqnarray}
\qed

\subsection{Stein's equation for asymmetric $\alpha$-stable distributions}
Let $h:\R\to\R$ be a Lispchitz function and let $N\sim N(0,1)$.
It is well-known and easy to prove (see, e.g., \cite[Prop. 3.5.1]{bluebook}) that the function
\[
f_h(x)=-\int_0^\infty \frac{e^{-t}}{\sqrt{1-e^{-2t}}}\,\E\left[
h(e^{-t}x+\sqrt{1-e^{-2t}}N)N
\right]dt
\]
is $C^1$ and
satisfies the Stein's equation associated with the standard gaussian distribution, namely $f'_h(x)-xf_h(x)=h(x)-\E[h(N)]$ for all $x\in\R$.

In this section, we introduce a function $\phi_h$ that satisfies an analogous property, but for the asymmetric $\alpha$-stable distribution (with $\alpha\in(1,2)$) instead of the gaussian one.
Because we want to keep our approach as elementary as possible, our proof
of Lemma \ref{lm23} is done `by hands', without relying on specific tools and results from the theory of semigroups.

\begin{lemma}\label{lm23}
Fix $\alpha\in(1,2)$ and $\beta\in[-1,1]$, and let  $h:\R\to\R$ be a Lipschitz function. Set
\begin{equation}\label{phih}
\phi_h(x) = -\int_0^\infty \int_\R p_{(1-e^{-t})^{\frac1\alpha},\beta}(y-e^{-\frac{t}{\alpha}}x)(h(y)-\mu(h))dydt,
\end{equation}
where $\mu$ is the distribution of $S_\alpha(1,\beta)$ (that is, $\mu(h)=\int_\R h(x)p_{1,\beta}(x)dx$).
Then
\begin{equation}\label{stablestein}
(\mathcal{L}^{\alpha,\beta} \phi_h)(x)-\frac{1}{\alpha}x\phi'_h(x) = h(x) - \mu(h).
\end{equation}
\end{lemma}
\noindent
{\it Proof}.
Since the identity (\ref{stablestein}) is linear with respect to $h$, it is enough (by approximation) to consider the case where $h=h_\lambda$ satisfies $h_\lambda'(x)=e^{i\lambda x}$ for some $\lambda\in\R\setminus\{0\}$.
First, one has
\begin{equation}
\phi'_{h_\lambda}(x)
=-
\int_0^\infty e^{-\frac{t}{\alpha}}\,e^{i\lambda e^{-\frac{t}{\alpha}}x}\,e^{-|\lambda|^\alpha(1-e^{-t})(1-i\,\beta\,{\rm sign}(\lambda)\,\tan\frac{\pi\alpha}2)} dt.
\label{phi'}
\end{equation}
Using Proposition \ref{properties} and then identity (\ref{def-den}) we deduce
\begin{eqnarray*}
(\mathcal{L}^{\alpha,\beta} \phi_{h_\lambda})(x)&=&\frac{d_\alpha}{\alpha}\int_0^\infty \frac{(1+\beta)\big(\phi'_{h_\lambda}(x+u)-\phi'_{h_\lambda}(x)\big)-(1-\beta)\big(\phi'_{h_\lambda}(x-u)-\phi'_{h_\lambda}(x)\big)}{2u^\alpha}du\\
&=&-\frac{d_\alpha}{\alpha}
\int_0^\infty
\frac{du}{2u^\alpha}
\int_0^\infty dt\,e^{-\frac{t}{\alpha}}e^{i\lambda e^{-\frac{t}{\alpha}}x}
e^{-|\lambda|^\alpha(1-e^{-t})(1-i\,\beta\,{\rm sign}(\lambda)\,\tan\frac{\pi\alpha}2)}\\
&&\hskip4cm\times
\left(
(1+\beta)(e^{i\lambda e^{-\frac{t}{\alpha}}u}-1)-(1-\beta)(
e^{-i\lambda e^{-\frac{t}{\alpha}}u}-1)
\right).
\end{eqnarray*}
Recall (see, e.g., \cite[identity (14.18)]{sato})
that $\int_0^\infty
(e^{ir}-1)r^{-1-\gamma}dr=\Gamma(-\gamma)e^{-\frac{i\pi\gamma}2}$
 for any $\gamma\in(0,1)$.
 Setting $\gamma=\alpha-1\in(0,1)$
 and doing the change of variable $v=|\lambda| u$ yield
 \[
 \int_0^\infty
(e^{i\lambda u}-1)u^{-\alpha}du
=|\lambda|^{\alpha -1}\int_0^\infty \big(e^{i\,{\rm sign}(\lambda)v}-1\big)v^{-\alpha}dv
=|\lambda|^{\alpha-1}\Gamma(1-\alpha)
e^{-i\,{\rm sign}(\lambda)\,\frac{\pi(\alpha-1)}{2}}.
\]
On the other hand,
\[
\frac1{d_\alpha}=
\int_0^\infty \frac{1-\cos y}{y^{1+\alpha}}dy=\frac1\alpha
\int_0^\infty \frac{\sin v}{v^\alpha} dx=\frac1\alpha\,\Gamma(1-\alpha)\cos\frac{\alpha\pi}{2},
\]
the last equality being obtained by applying \cite[identity 3.764]{GR}.
We deduce
\begin{eqnarray*}
&& \frac{d_\alpha}{2\alpha}\int_0^\infty
\left(
(1+\beta)(e^{i\lambda u}-1)-(1-\beta)(
e^{-i\lambda u}-1)
\right)u^{-\alpha}du\\
&&\hskip3cm =\frac{|\lambda|^{\alpha-1}}{2\cos\frac{\pi\alpha}{2}}
\left(
(1+\beta)e^{-i\,{\rm sign}(\lambda)\,\frac{\pi(\alpha-1)}{2}} - (1-\beta)
e^{i\,{\rm sign}(\lambda)\,\frac{\pi(\alpha-1)}{2} }\right)\\
&&\hskip3cm =i\,{\rm sign}(\lambda)|\lambda|^{\alpha-1}\left(
1-i\,\beta\,{\rm sign}(\lambda)\,\tan\frac{\pi\alpha}2\right),
\end{eqnarray*}
implying in turn
\begin{eqnarray*}
&&(\mathcal{L}^{\alpha,\beta} \phi_{h_\lambda})(x)\\
&&\hskip1cm=-i\,\frac{|\lambda|^{\alpha}}{\lambda}\,\left(
1-i\,\beta\,{\rm sign}(\lambda)\,\tan\frac{\pi\alpha}2\right)
\,
\int_0^\infty e^{i\lambda e^{-\frac{t}{\alpha}}x}
e^{-|\lambda|^\alpha(1-e^{-t})(1-i\,\beta\,{\rm sign}(\lambda)\,\tan\frac{\pi\alpha}2)}\,e^{-t}dt.
\end{eqnarray*}
Finally, integrating by parts with $u(t)=e^{i\lambda e^{-\frac{t}{\alpha}}x}$ and $v'(t)=e^{-t}\,e^{-|\lambda|^\alpha(1-e^{-t})(1-i\,\beta\,{\rm sign}(\lambda)\,\tan\frac{\pi\alpha}2)}$ yields, using also (\ref{phi'}),
\begin{eqnarray*}
(\mathcal{L}^{\alpha,\beta} \phi_{h_\lambda})(x)
&=&\frac{1}{i\lambda}\big(e^{i\lambda x}-e^{-|\lambda|^\alpha(1-i\,\beta\,{\rm sign}(\lambda)\,\tan\frac{\pi\alpha}2)}\big)
+\frac{x}{\alpha}\phi'_{h_\lambda}(x).
\end{eqnarray*}
Since $h_\lambda(x)-\mu(h_\lambda)=\frac{1}{i\lambda}\big(e^{i\lambda x}-e^{-|\lambda|^\alpha(1-i\,\beta\,{\rm sign}(\lambda)\,\tan\frac{\pi\alpha}2)})$,
the desired conclusion (\ref{stablestein}) follows.
\qed

\subsection{Proof of Theorem \ref{thm1}}
We are now ready to proceed with the proof of Theorem \ref{thm1}.
Recall $p_{\sigma,\beta}$ is the density of $S_{\alpha}(\sigma,\beta)$.
By scaling, we observe that $p_{\sigma,\beta}(x)=\frac1\sigma p_{1,\beta}(\frac{x}{\sigma})$ for all $\sigma>0$ and $x\in\R$, implying in turn that
$p'_{\sigma,\beta}(x)=\frac1{\sigma^2} p'_{1,\beta}(\frac{x}{\sigma})$
for all $\sigma>0$ and $x\in\R$.

On the other hand, by Fourier inversion we have
\[
p_{1,\beta}(x)=\frac1{2\pi}\int_\R e^{-i\lambda x -
|\lambda|^\alpha(1-i\,\beta\,{\rm sign}(\lambda)\,\tan\frac{\pi\alpha}2)}d\lambda,
\]
implying in turn that
\[
p'_{1,\beta}(x)=\frac{-i}{2\pi}\int_\R \lambda e^{-i\lambda x -
|\lambda|^\alpha(1-i\,\beta\,{\rm sign}(\lambda)\,\tan\frac{\pi\alpha}2)}d\lambda.
\]
As a result,
$
\|p'_{1,\beta}\|_\infty <\infty.
$
Using two integrations by parts, one proves that
$x\mapsto x^2p'_{1,\beta}(x)$ is bounded too; these two facts together
implies that $\int_\R |p'_{1,\beta}(x)|dx<\infty$.

Now, fix a Lipschitz function $h:\R\to\R$ and recall $\phi_h$ from (\ref{phih}).
We observe that
\begin{eqnarray*}
\big|\phi_h'(x)\big| &=& \left|\int_0^\infty dt\, e^{-\frac{t}{\alpha}}\int_\R dy\,p'_{(1-e^{-t})^\frac1\alpha,\beta}(y-e^{-\frac{t}{\alpha}}x)(h(y)-\mu(h))\right|\\
&=& \left| \int_0^\infty dt\,e^{-\frac{t}{\alpha}}\int_\R dy\,p_{(1-e^{-t})^\frac1\alpha,\beta}(y-e^{-\frac{t}{\alpha}}x)h'(y)\right|\\
&\leq&\|h'\|_\infty \int_0^\infty e^{-\frac{t}{\alpha}}dt \int_\R p_{(1-e^{-t})^\frac1\alpha,\beta}(y)dy=
\alpha \|h'\|_\infty,
\end{eqnarray*}
whereas
\begin{eqnarray*}
\big|\phi_h''(x)\big| &=& \left| \int_0^\infty dt\,e^{-\frac{2t}{\alpha}}\int_\R dy\,p'_{(1-e^{-t})^\frac1\alpha,\beta}(y-e^{-\frac{t}{\alpha}}x)h'(y)\right|\\
&\leq&\|h'\|_\infty \int_0^\infty e^{-\frac{2t}{\alpha}}dt \int_\R |p'_{(1-e^{-t})^\frac1\alpha,\beta}(y)|dy\\
&=&\|h'\|_\infty \int_0^1 u^{\frac{2}{\alpha}-1}(1-u)^{-\frac1\alpha}du \int_\R |p'_{1,\beta}(y)|dy=\eta_{\alpha,\beta} \,\|h'\|_\infty.
\end{eqnarray*}

To conclude, it now suffices to consider (\ref{stablestein}) with $x=X$,
to take expectation in both sides, and to use the two previous bounds.
\qed

\section{Proof of Theorem \ref{thm2} and Theorem \ref{thm3}}

In order to prove Theorem \ref{thm2} and Theorem \ref{thm3}, we extend the celebrated Stein's leave-one-out approach classically used in
the context of normal approximation (see, e.g., \cite[pages 5-6]{CGS}).

\subsection{Taylor-like extension}

We shall make use of the following lemmas.
\begin{lemme}\label{lm41}
Let $X$ have a distribution of the form (\ref{px}) with $\epsilon(x)$ satisfying (\ref{bound}) and
$Y$ be two independent integrable random variables.
For any $0<a<(2A)^{-\frac{1}{\alpha}}\wedge1$
and any $\phi:\R\to\R$ such that $\|\phi'\|_\infty,\|\phi''\|_\infty<\infty$, denote
\begin{align*}
T=&\left|
\E[X\phi'(Y+aX)] -\E[X]\E[\phi'(Y)]-\frac{2A\alpha^{2}}{d_\alpha}\,a^{\alpha-1}\E\left[\big(\mathcal{L}^{\alpha,\beta}\phi\big)(Y)\right]
\right|.
\end{align*}
Then:
\begin{enumerate}
\item[i)] When $\gamma\in(2-\alpha,\infty),$ we have
\begin{align*}
T\leq2(2A)^{\frac{2}{\alpha}}\Big[\frac{2}{2-\alpha}+\frac{2K}{\alpha+\gamma-2}(2A)^{\frac{-\alpha-\gamma}{\alpha}}\Big]\|\phi''\|_{\infty}a.
\end{align*}
\item[ii)] When $\gamma=2-\alpha,$ we have
\begin{align*}
T\leq\frac{2\alpha}{2-\alpha}(2A)^{\frac{2}{\alpha}}\|\phi''\|_{\infty}a+\Big[\Big(2(2A)^{\frac{2}{\alpha}}+\frac{8K}{\alpha-1}\Big)\|\phi''\|_{\infty}+\frac{8\alpha(A+K)-4K}{\alpha-1}\|\phi'\|_{\infty}\Big]a\,|\log a|.
\end{align*}
\item[iii)] When $\gamma\in(0,2-\alpha),$ we have
\begin{align*}
T\leq\Big[\Big(\frac{4(2A)^{\frac{2}{\alpha}}}{2-\alpha}+\frac{8K}{2-\alpha-\gamma}\Big)
\|\phi''\|_{\infty}+\frac{8\alpha(A+K)-4K}{\alpha-1}\|\phi'\|_{\infty}\Big]a^{\frac{1-\alpha}{\gamma-1}}.
\end{align*}
\item[iv)] When $\gamma=0,$ we have
\begin{align*}
T\leq&\frac{2\alpha(2A)^{\frac{2}{\alpha}}}{2-\alpha}\|\phi''\|_{\infty}a+4\|\phi''\|_{\infty}a\int_{-a^{-1}}^{a^{-1}}\frac{|\epsilon(x)|}{|x|^{\alpha-1}}dx\\
&+\Big[\Big(\frac{8}{2-\alpha}+2(2A)^{\frac{2}{\alpha}}\Big)\|\phi''\|_{\infty}+\frac{8\alpha(A+K)-4K}{\alpha-1}\|\phi'\|_{\infty}\Big]a^{\alpha-1}\big(\sup_{|x|\geq a^{-1}}|\epsilon(x)|\big)^{\alpha-1}.
\end{align*}
\end{enumerate}
\end{lemme}
\noindent
{\it Proof}. We can write, using (\ref{crucial}),
\begin{align*}
&\frac{2A\alpha^{2}}{d_{\alpha}}a^{\alpha-1}\mathbb{E}\left[\big(\mathcal{L}^{\alpha,\beta}\phi\big)(Y)\right]\\
=&2A\alpha\mathbb{E}\Big[\int_{-\infty}^{\infty}u\big[\phi'(Y+au)-\phi'(Y)\big]\frac{(1+\beta){\bf 1}_{(0,\infty)}(u)+(1-\beta){\bf 1}_{(-\infty,0)}(u)}
{2|u|^{\alpha+1}}du\Big]\\
=&\mathbb{E}\Big[\int_{|u|\geq(2A)^{\frac{1}{\alpha}}}u\big[\phi'(Y+au)-\phi'(Y)\big]\frac{A\alpha\big[(1+\beta){\bf 1}_{(0,\infty)}(u)+(1-\beta){\bf 1}_{(-\infty,0)}(u)\big]}
{|u|^{\alpha+1}}du\Big]+\mathcal{R},
\end{align*}
where
\begin{align}\label{R}
\mathcal{R}=\mathbb{E}\Big[\int_{-(2A)^{\frac{1}{\alpha}}}^{(2A)^{\frac{1}{\alpha}}}u\big[\phi'(Y+au)-\phi'(Y)\big]\frac{A\alpha\big[(1+\beta){\bf 1}_{(0,\infty)}(u)+(1-\beta){\bf 1}_{(-\infty,0)}(u)\big]}
{|u|^{\alpha+1}}du\Big].
\end{align}
Since $\int_{|u|\geq(2A)^{\frac{1}{\alpha}}}\frac{2A\alpha\big[(1+\beta){\bf 1}_{(0,\infty)}(u)+(1-\beta){\bf 1}_{(-\infty,0)}(u)\big]}
{2|u|^{\alpha+1}}du=1,$ we can consider a random variable $\tilde{X}$ which is independent of $Y$ and satisfies
\begin{align}\label{F}
\mathbb{P}(\tilde{X}>x)=\frac{A(1+\beta)}{|x|^{\alpha}},\quad x\geq(2A)^{\frac{1}{\alpha}},\qquad \mathbb{P}(\tilde{X}\leq x)=\frac{A(1-\beta)}{|x|^{\alpha}},\quad x\leq-(2A)^{\frac{1}{\alpha}}.
\end{align}
We then have
\begin{align*}
\frac{2A\alpha^{2}}{d_{\alpha}}a^{\alpha-1}\mathbb{E}\left[\big(\mathcal{L}^{\alpha,\beta}\phi\big)(Y)\right]=\mathbb{E}\big[\tilde{X}\phi'(Y+a\tilde{X})\big]-\mathbb{E}[\tilde{X}]\mathbb{E}\big[\phi'(Y)\big]
+\mathcal{R}.
\end{align*}
As a result, denoting by $F_{\tilde{X}}$ the distribution function of $\tilde{X}$,
we have
\begin{align}\label{result}
&\Big|\E[X\phi'(Y+aX)] -\E[X]\E[\phi'(Y)]-\frac{2A\alpha^{2}}{d_\alpha}\,a^{\alpha-1}\E\left[\big(\mathcal{L}^{\alpha,\beta}\phi\big)(Y)\right]\Big|\nonumber\\
\leq&\mathbb{E}\Big|\int_{-\infty}^{\infty}x\big[\phi'(Y+ax)-\phi'(Y)\big]
d\big(F_{X}(x)-F_{\tilde{X}}(x)\big)\Big|+|\mathcal{R}|
\end{align}
and
\begin{align*}
F_{X}(x)-F_{\tilde{X}}(x)=&\big(\frac{1}{2}-\frac{A+\epsilon(x)}{|x|^{\alpha}}\big)(1+\beta){\bf 1}_{(0,(2A)^{\frac{1}{\alpha}})}(x)-\frac{\epsilon(x)}{|x|^{\alpha}}(1+\beta){\bf 1}_{((2A)^{\frac{1}{\alpha}},\infty)}(x)\\
&+\big(\frac{A+\epsilon(x)}{|x|^{\alpha}}-\frac{1}{2}\big)(1-\beta){\bf 1}_{(-(2A)^{\frac{1}{\alpha}},0)}(x)+\frac{\epsilon(x)}{|x|^{\alpha}}(1-\beta){\bf 1}_{(-\infty,-(2A)^{\frac{1}{\alpha}})}(x).
\end{align*}
It is easy to verify
\begin{align}\label{R1}
|\mathcal{R}|\leq A\alpha\|\phi''\|_{\infty}a\int_{-(2A)^{\frac{1}{\alpha}}}^{(2A)^{\frac{1}{\alpha}}}\frac{(1+\beta){\bf 1}_{(0,\infty)}(u)+(1-\beta){\bf 1}_{(-\infty,0)}(u)}
{|u|^{\alpha-1}}du=\frac{2\alpha}{2-\alpha}(2A)^{\frac{2}{\alpha}}\|\phi''\|_{\infty}a.
\end{align}
Now, let us deal with the first term of (\ref{result}). Recall our assumption (\ref{bound}). We split into two different cases, according to the place of
$\gamma$ with respect to $2-\alpha$.\\
1. Assume $\gamma>2-\alpha.$ We have, by integrating by parts
\begin{align*}
&\mathbb{E}\Big|\int_{-\infty}^{\infty}x\big[\phi'(Y+ax)-\phi'(Y)\big]
d\big(F_{X}(x)-F_{\tilde{X}}(x)\big)\Big|\\
=&\mathbb{E}\Big|\int_{-\infty}^{\infty}\big(F_{X}(x)-F_{\tilde{X}}(x)\big)\big[ax\phi''(Y+ax)+\phi'(Y+ax)-\phi'(Y)\big]
dx\Big|\\
\leq&2\|\phi''\|_{\infty}a\Big[\int_{-(2A)^{\frac{1}{\alpha}}}^{(2A)^{\frac{1}{\alpha}}}|x|dx+\int_{(2A)^{\frac{1}{\alpha}}}^{\infty}\frac{|\epsilon(x)|}{x^{\alpha-1}}dx+
\int_{-\infty}^{-(2A)^{\frac{1}{\alpha}}}\frac{|\epsilon(x)|}{|x|^{\alpha-1}}dx\Big]\\
\leq&2\|\phi''\|_{\infty}a\Big[(2A)^{\frac{2}{\alpha}}+2\int_{(2A)^{\frac{1}{\alpha}}}^{\infty}\frac{K}{x^{\alpha+\gamma-1}}dx\Big]=2(2A)^{\frac{2}{\alpha}}
\big[1+\frac{2K}{\alpha+\gamma-2}(2A)^{\frac{-\alpha-\gamma}{\alpha}}\big]\|\phi''\|_{\infty}a.
\end{align*}
2. Assume now $0\leq\gamma\leq2-\alpha$. Choose a number $N>a^{-1}.$ One has by \cite[Lemma 2.8]{lihu} and using that $|\epsilon(x)|\leq K$ for $|x|>N$,
\begin{align*}
&\mathbb{E}\Big|\int_{|x|>N}x\big[\phi'(Y+ax)-\phi'(Y)\big]
d\big(F_{X}(x)-F_{\tilde{X}}(x)\big)\Big|\\
\leq&2\|\phi'\|_{\infty}\Big[\int_{|x|>N}|x|dF_{X}(x)+\int_{|x|>N}|x|dF_{\tilde{X}}(x)\Big]\\
=&2\|\phi'\|_{\infty}\mathbb{E}\big[|X|{\bf 1}_{(N,\infty)}(|X|)+|\tilde{X}|{\bf 1}_{(N,\infty)}(|\tilde{X}|)\big]\leq\frac{4(2A+K)\alpha}{\alpha-1}\|\phi'\|_{\infty}N^{1-\alpha}.
\end{align*}
On the other hand, by integrating by parts
\begin{align*}
&\mathbb{E}\Big|\int_{-N}^{N}x\big[\phi'(Y+ax)-\phi'(Y)\big]
d\big(F_{X}(x)-F_{\tilde{X}}(x)\big)\Big|\\
\leq&4K\|\phi'\|_{\infty}N^{1-\alpha}+2\|\phi''\|_{\infty}a\int_{-N}^{N}\big|xF_{X}(x)-xF_{\tilde{X}}(x)\big|dx\\
\leq&4K\|\phi'\|_{\infty}N^{1-\alpha}+2\|\phi''\|_{\infty}a\big[(2A)^{\frac{2}{\alpha}}+2\int_{(2A)^{\frac{1}{\alpha}}}^{N}\frac{|\epsilon(x)|}{x^{\alpha-1}}dx
+2\int^{-(2A)^{\frac{1}{\alpha}}}_{-N}\frac{|\epsilon(x)|}{|x|^{\alpha-1}}dx\big].
\end{align*}
If $0<\gamma\leq2-\alpha,$ we have
\begin{align*}
\int_{(2A)^{\frac{1}{\alpha}}}^{N}\frac{|\epsilon(x)|}{x^{\alpha-1}}dx\leq\int_{(2A)^{\frac{1}{\alpha}}}^{N}\frac{K}{x^{\alpha+\gamma-1}}dx\leq
\begin{cases}
K\log N,\quad &\gamma=2-\alpha,\\
\frac{K}{2-\alpha-\gamma}N^{2-\alpha-\gamma},\quad &\gamma\in(0,2-\alpha).
\end{cases}
\end{align*}
If $\gamma=0,$ we have
\begin{align*}
\int_{(2A)^{\frac{1}{\alpha}}}^{N}\!\frac{|\epsilon(x)|}{x^{\alpha-1}}dx=\!\int_{(2A)^{\frac{1}{\alpha}}}^{a^{-1}}\!\frac{|\epsilon(x)|}{x^{\alpha-1}}dx+\int_{a^{-1}}^{N}\!\frac{|\epsilon(x)|}{x^{\alpha-1}}dx
\leq\!\int_{0}^{a^{-1}}\!\frac{|\epsilon(x)|}{x^{\alpha-1}}dx+\!\frac{1}{2-\alpha}\sup_{x>a^{-1}}|\epsilon(x)|N^{2-\alpha}.
\end{align*}
Since similar bounds hold true for $\int^{-(2A)^{\frac{1}{\alpha}}}_{-N}\frac{|\epsilon(x)|}{|x|^{\alpha-1}}dx,$ we can consider
\begin{align*}
N^{1-\alpha}=
\begin{cases}
aN^{2-\alpha-\gamma},\quad &\gamma\in(0,2-\alpha]\\
a\sup_{|x|>a^{-1}}|\epsilon(x)|N^{2-\alpha},\quad &\gamma=0,
\end{cases}
\end{align*}
which implies
\begin{align*}
N=
\begin{cases}
a^{\frac{1}{\gamma-1}},\quad &\gamma\in(0,2-\alpha],\\
a^{-1}\big(\sup_{|x|>a^{-1}}|\epsilon(x)|)^{-1},\quad &\gamma=0.
\end{cases}
\end{align*}
The desired conclusion follows.
\qed

\begin{lemme}\label{lm42}
Keep the same notation and assumptions as in Lemma \ref{lm41}. In addition, further assume $\frac{\epsilon(x)}{x^{\alpha}}{\bf 1}_{(0,\infty)}(x)$ and $\frac{\epsilon(x)}{|x|^{\alpha}}{\bf 1}_{(-\infty,0)}(x)$ are ultimately monotone. We have, for any $a>0$ such that $\frac{\epsilon(x)}{x^{\alpha}}{\bf 1}_{(0,\infty)}(x)$ and $\frac{\epsilon(x)}{|x|^{\alpha}}{\bf 1}_{(-\infty,0)}(x)$ are monotone for $|x|>a^{-1}$ and any $\phi:\R\to\R$ such that $\|\phi'\|_\infty,\|\phi''\|_\infty<\infty$:
\begin{eqnarray*}
&&\left|
\E[X\phi'(Y+aX)] -\E[X]\E[\phi'(Y)]-\frac{2A\alpha^{2}}{d_\alpha}\,a^{\alpha-1}\E\left[\big(\mathcal{L}^{\alpha,\beta}\phi\big)(Y)\right]
\right|\\
&&\hskip2cm\leq2(2A)^{\frac{1}{\alpha}}a+
\frac{(16\alpha-1)\|\phi'\|_\infty}{\alpha-1}\,a^{\alpha-1}\sup_{|x|\geq a^{-1}}|\epsilon(x)|+
2\|\phi''\|_\infty\,a\int_{-1/a}^{1/a}\frac{|\epsilon(x)|}{|x|^{\alpha-1}}dx.
\end{eqnarray*}
\end{lemme}
\noindent
{\it Proof}.
By (\ref{result}), we have
\begin{align*}
&\Big|\E[X\phi'(Y+aX)] -\E[X]\E[\phi'(Y)]-\frac{2A\alpha^{2}}{d_\alpha}\,a^{\alpha-1}\E\left[\big(\mathcal{L}^{\alpha,\beta}\phi\big)(Y)\right]\Big|\nonumber\\
\leq&\mathbb{E}\Big|\int_{-\infty}^{\infty}x\big[\phi'(Y+ax)-\phi'(Y)\big]
d\big(F_{X}(x)-F_{\tilde{X}}(x)\big)\Big|+|\mathcal{R}|,
\end{align*}
where $F_{\tilde{X}}(x)$ and $\mathcal{R}$ are defined in (\ref{F}) and (\ref{R}), respectively. By (\ref{R1}), we know
\begin{align*}
|\mathcal{R}|\leq\frac{2\alpha}{2-\alpha}(2A)^{\frac{2}{\alpha}}\|\phi''\|_{\infty}a.
\end{align*}
For the first term, one has, by integrating by parts
\begin{align*}
&\mathbb{E}\Big|\int_{a^{-1}}^{\infty}x\big[\phi'(Y+ax)-\phi'(Y)\big]
d\big(F_{X}(x)-F_{\tilde{X}}(x)\big)\Big|\\
\leq&2\|\phi'\|_{\infty}\int_{a^{-1}}^{\infty}x\big|d\big(F_{X}(x)-F_{\tilde{X}}(x)\big)\big|\\
\leq&4\|\phi'\|_{\infty}\Big(a^{\alpha-1}|\epsilon(a^{-1})|+\int_{a^{-1}}^{\infty}\frac{|\epsilon(x)|}{x^{\alpha}}dx\Big)\leq\frac{4\alpha}{\alpha-1}\|\phi'\|_{\infty}a^{\alpha-1}\sup_{x\geq a^{-1}}|\epsilon(x)|.
\end{align*}
In the same way
\begin{align*}
\mathbb{E}\Big|\int^{-a^{-1}}_{-\infty}x\big[\phi'(Y+ax)-\phi'(Y)\big]
d\big(F_{X}(x)-F_{\tilde{X}}(x)\big)\Big|\leq\frac{4\alpha}{\alpha-1}\|\phi'\|_{\infty}a^{\alpha-1}\sup_{x\leq-a^{-1}}|\epsilon(x)|.
\end{align*}
On the other hand,
\begin{align*}
&\mathbb{E}\Big|\int_{-a^{-1}}^{a^{-1}}x\big[\phi'(Y+ax)-\phi'(Y)\big]
d\big(F_{X}(x)-F_{\tilde{X}}(x)\big)\Big|\\
\leq&4\|\phi'\|_{\infty}a^{\alpha-1}|\epsilon(a^{-1})|+4\|\phi'\|_{\infty}a^{\alpha-1}|\epsilon(-a^{-1})|+2(2A)^{\frac{1}{\alpha}}a+2\|\phi''\|_{\infty}a\int_{-a^{-1}}^{a^{-1}}\frac{|\epsilon(x)|}{|x|^{\alpha-1}}dx,
\end{align*}
and the desired conclusion follows.\qed

\subsection{Leave-one out method and proof of Theorem \ref{thm2} and \ref{thm3}}\label{loo}
Recall the notation introduced in Theorem \ref{thm2}.
We have $\sigma=\left(\frac{2A\alpha}{d_\alpha}\right)^\frac1\alpha$
and $S_{n,i}=S_n-\frac{n^{-\frac{1}{\alpha}}}{\sigma}(X_i-\E[X_i])$. By observing that $S_{n,i}$ and $X_i$ are independent,
one can write:
\[
\left|\E[S_n\phi'(S_n)] -\alpha\,\E[(\mathcal{L}^{\alpha,\beta}\phi)(S_n)]\right|\leq I+II+III+IV,
\]
where
\begin{eqnarray*}
I&=&\frac\alpha{n}\sum_{i=1}^n
\Big|
\E[(\mathcal{L}^{\alpha,\beta}\phi)(S_{n,i})]-
\E[(\mathcal{L}^{\alpha,\beta}\phi)(S_{n})]\Big|
\\
II&=&\frac{n^{-\frac{1}{\alpha}}}{\sigma}
\sum_{i=1}^n
\left|
\E\left[
X_i\,\phi'\left(
S_{n,i}+\frac{n^{-\frac{1}{\alpha}}}{\sigma}X_i
\right)
\right]\!\!-\!\E[X_i]\E[\phi'(
S_{n,i})]
-\frac{2A\alpha^{2}}{d_\alpha}\left(
\frac{n^{-\frac1\alpha}}{\sigma}
\right)^{\alpha-1}
\!\!\!\!\!\!\!\!\E[(\mathcal{L}^{\alpha,\beta}\phi)(S_{n,i})]
\right|\\
III&=&\frac{n^{-\frac{1}{\alpha}}}{\sigma}
\sum_{i=1}^n
\left|
\E\left[
X_i\left(
\phi'\left(
S_{n,i}+\frac{n^{-\frac{1}{\alpha}}}{\sigma}(X_i-\E[X_i])
\right)
-\phi'\left(
S_{n,i}+\frac{n^{-\frac{1}{\alpha}}}{\sigma}X_i
\right)
\right)
\right]
\right|\\
IV&=&\frac{n^{-\frac{1}{\alpha}}}{\sigma}
\sum_{i=1}^n
\big|\E[X_i]\big|\,\left|\E\left[
\phi'(
S_{n,i})
-\phi'\left(
S_{n,i}+\frac{n^{-\frac{1}{\alpha}}}{\sigma}(X_i-\E[X_i])
\right)
\right]\right|.
\end{eqnarray*}
We have, thanks to (\ref{crucialbis}):
\begin{eqnarray*}
I
&\leq&
\frac{4d_{\alpha}}{(2-\alpha)(\alpha-1)}\|\phi''\|_\infty\,\frac{\E[|X_1-\E[X_1]|^{2-\alpha}]}{\sigma^{2-\alpha}}\,n^{1-\frac2\alpha}.
\end{eqnarray*}
Using Lemma \ref{lm41},
\begin{enumerate}
\item[i)] When $\gamma\in(2-\alpha,\infty),$ we have
\begin{eqnarray*}
II&\leq&\frac{2(2A)^{\frac{2}{\alpha}}}{\sigma^{2}}\Big[\frac{2}{2-\alpha}+\frac{2K}{\alpha+\gamma-2}(2A)^{\frac{-\alpha-\gamma}{\alpha}}\Big]\|\phi''\|_{\infty}n^{\frac{\alpha-2}{\alpha}}.
\end{eqnarray*}
\item[ii)] When $\gamma=2-\alpha,$ we have
\begin{eqnarray*}
II&\leq&\frac{1}{\sigma^{2}}\Big[\Big(\frac{4(2A)^{\frac{2}{\alpha}}}{2-\alpha}+\frac{8K}{\alpha-1}\Big)\|\phi''\|_{\infty}+\frac{8\alpha(A+K)-4K}{\alpha-1}\|\phi'\|_{\infty}\Big]n^{\frac{\alpha-2}{\alpha}}|\log(\sigma n^{\frac{1}{\alpha}})|.
\end{eqnarray*}
\item[iii)] When $\gamma\in(0,2-\alpha),$ we have
\begin{eqnarray*}
II&\leq&\sigma^{\frac{\alpha-\gamma}{\gamma-1}}\Big[\Big(\frac{4(2A)^{\frac{2}{\alpha}}}{2-\alpha}+\frac{8K}{2-\alpha-\gamma}\Big)
\|\phi''\|_{\infty}+\frac{8\alpha(A+K)-4K}{\alpha-1}\|\phi'\|_{\infty}\Big]n^{-\frac{(\alpha-1)\gamma}{\alpha(1-\gamma)}}.
\end{eqnarray*}
\item[iv)] When $\gamma=0,$ we have
\begin{eqnarray*}
II&\leq&\frac{2\alpha(2A)^{\frac{2}{\alpha}}}{(2-\alpha)\sigma^{2}}\|\phi''\|_{\infty}n^{\frac{\alpha-2}{\alpha}}+\frac{4\|\phi''\|_{\infty}}{\sigma^{2}}n^{\frac{\alpha-2}{\alpha}}\int_{-\sigma n^{\frac{1}{\alpha}}}^{\sigma n^{\frac{1}{\alpha}}}\frac{|\epsilon(x)|}{|x|^{\alpha-1}}dx\\
&&+\sigma^{-\alpha}\Big[\Big(\frac{8}{2-\alpha}+2(2A)^{\frac{2}{\alpha}}\Big)\|\phi''\|_{\infty}+\frac{8\alpha(A+K)-4K}{\alpha-1}\|\phi'\|_{\infty}\Big]\big(\sup_{|x|\geq \sigma n^{\frac{1}{\alpha}}}|\epsilon(x)|\big)^{\alpha-1}.
\end{eqnarray*}
\end{enumerate}
Finally, by the mean value theorem:
\begin{eqnarray*}
III+IV&\leq&\frac{3\|\phi''\|_\infty}{\sigma^2}\E[|X_1|]|\E[X_1]|\,n^{1-\frac2\alpha}.
\end{eqnarray*}
Plugging this into (\ref{steinstable}) gives the desired conclusion of Theorem \ref{thm2}.\qed

\begin{proof} [{\bf Proof of Theorem \ref{thm3}}]
It suffices to bound the $II$ in the proof of Theorem \ref{thm2}. Using Lemma \ref{lm42}, we have
\begin{eqnarray*}
II&\leq&\frac{2(2A)^{\frac{1}{\alpha}}}{\sigma^{2}}n^{1-\frac{2}{\alpha}}+\frac{(16\alpha-1)\|\phi'\|_\infty}{(\alpha-1)\sigma^\alpha}\,
\sup_{|x|\geq \sigma\,n^{\frac1\alpha}}|\epsilon(x)|+
\frac{2\|\phi''\|_\infty}{\sigma^2}\,n^{1-\frac{2}{\alpha}}\int_{-\sigma\,n^{\frac1\alpha}}^{\sigma\,n^{\frac1\alpha}}
\frac{|\epsilon(x)|}{|x|^{\alpha-1}}dx,
\end{eqnarray*}
from which the desired conclusion follows.
\end{proof}
\section{A more difficult example: proof of (\ref{ex4})}
In this section, we prove the estimate (\ref{ex4}).
Consider independent copies $X_1,\ldots,X_n$ of a random variable with density
$
p_X(x)=\frac{\alpha^2e^\alpha}{2(1+\alpha)}\,\frac{\log |x|}{|x|^{\alpha+1}}{\bf 1}_{[e,\infty)}(|x|)
$
and define the sequence $(\gamma_n)_{n\geq 1}$ implicitly by $\gamma_n=\big(n\log \gamma_n\big)^\frac1\alpha$.
Observe that $X_1$ is integrable and centered.
We set $\sigma =\left(\frac{\alpha^2e^\alpha}{(1+\alpha)d_\alpha}\right)^{\frac1\alpha}$, $\widetilde{X}_i=\frac{n^\frac1\alpha}{\gamma_n}X_i$, $\widetilde{S}_n=\frac1\sigma\,n^{-\frac1\alpha}(\widetilde{X}_1+\ldots+\widetilde{X}_n)$, and
$\widetilde{S}_{n,i}=\widetilde{S}_n-\frac{n^{-\frac1\alpha}}{\sigma}\widetilde{X}_i$.

To prove (\ref{ex4}), we shall use Theorem \ref{thm1} with $\beta=0$.
Let $\phi\in C^2(\R)$ be such that $\|\phi'\|_\infty\leq \alpha$ and
$\|\phi''\|_\infty\leq \eta_{\alpha,0}$.
We can write
\begin{eqnarray*}
&&\E[\widetilde{S}_n\phi'(\widetilde{S}_n)]-\alpha\,\E\big[
(\mathcal{L}^{\alpha,0}\phi)(\widetilde{S}_n)\big]\\
&=&\frac\alpha n \sum_{i=1}^n \left(
\E\big[
(\mathcal{L}^{\alpha,0}\phi)(\widetilde{S}_{n,i})\big]
-\E\big[
(\mathcal{L}^{\alpha,0}\phi)(\widetilde{S}_n)\big]
\right)\\
&&+\frac{n^{-\frac1\alpha}}{\sigma}\sum_{i=1}^n
\left(
\E\left[
\widetilde{X}_i\,\phi'(\widetilde{S}_{n,i}+\frac{n^{-\frac1\alpha}}{\sigma}\widetilde{X}_i)
\right]
-\frac{\alpha^3e^\alpha}{d_\alpha(1+\alpha)}\sigma^{1-\alpha}n^{\frac1\alpha-1}\,\E\big[
(\mathcal{L}^{\alpha,0}\phi)(\widetilde{S}_{n,i})\big]
\right)
\end{eqnarray*}

We have, using among other that $n\gamma_n^{-\alpha}=\frac1{\log \gamma_n}$,
\begin{eqnarray*}
&&\E\left[
\widetilde{X}_i\,\phi'(\widetilde{S}_{n,i}+\frac{n^{-\frac1\alpha}}{\sigma}\widetilde{X}_i)\right]\\
&=&\frac{\alpha^2e^\alpha}{2(1+\alpha)}\,\E\left[\int_\R \big(\phi'(\widetilde{S}_{n,i}+\frac{n^{-\frac1\alpha}}{\sigma}u)-\phi'(\widetilde{S}_{n,i})\big)
\frac{
u\log(
n^{-\frac1\alpha} \gamma_n\big|u|
)
}
{
|u|^{\alpha +1}\log \gamma_n
}\,
{\bf 1}_{[e,\infty)}(n^{-\frac1\alpha}\gamma_n\,|u|)\,
du\right]\\
&=&\frac{\alpha^2e^\alpha}{2(1+\alpha)}\,\E\left[\int_\R \big(\phi'(\widetilde{S}_{n,i}+\frac{n^{-\frac1\alpha}}{\sigma}u)-\phi'(\widetilde{S}_{n,i})\big)
\frac{
u
}
{
|u|^{\alpha +1}
}\,
{\bf 1}_{[e,\infty)}(n^{-\frac1\alpha}\gamma_n\,|u|)\,
du\right]\\
&&+\frac{\alpha^2e^\alpha}{2(1+\alpha)}\,\E\left[\int_\R \big(\phi'(\widetilde{S}_{n,i}+\frac{n^{-\frac1\alpha}}{\sigma}u)-\phi'(\widetilde{S}_{n,i})\big)
\frac{
u\log(
n^{-\frac1\alpha} \big|u|
)
}
{
|u|^{\alpha +1}\log \gamma_n
}\,
{\bf 1}_{[e,\infty)}(n^{-\frac1\alpha}\gamma_n\,|u|)\,
du\right].
\end{eqnarray*}
On the other hand, the identity (\ref{crucial}) with $a=\frac{n^{-\frac1\alpha}}{\sigma}$ yields
\begin{eqnarray*}
\frac{2\alpha}{d_\alpha}\,\sigma^{1-\alpha}n^{\frac1\alpha-1}(\mathcal{L}^{\alpha,0}\phi)(\widetilde{S}_{n,i})
&=&\int_\R
\big(\phi'(\widetilde{S}_{n,i}+\frac{n^{-\frac1\alpha}}{\sigma}u)-\phi'(\widetilde{S}_{n,i})\big)\frac{u}{|u|^{\alpha +1}}du.
\end{eqnarray*}
As a result,
\begin{eqnarray*}
&&\left|\E\left[
\widetilde{X}_i\,\phi'(\widetilde{S}_{n,i}+\frac{n^{-\frac1\alpha}}{\sigma}\widetilde{X}_i)\right]- \frac{\alpha^3e^\alpha}{d_\alpha(1+\alpha)}\,\sigma^{1-\alpha}\,n^{\frac1\alpha-1}\,\E\left[(\mathcal{L}^{\alpha,0}\phi)(\widetilde{S}_{n,i})\right]\right|\\
&=&\left|\frac{\alpha^2e^\alpha}{2(1+\alpha)}\,\E\left[\int_\R \big(\phi'(\widetilde{S}_{n,i}+\frac{n^{-\frac1\alpha}}{\sigma}u)-\phi'(\widetilde{S}_{n,i})\big)
\frac{
u
}
{
|u|^{\alpha +1}
}\,
{\bf 1}_{(0,e)}(n^{-\frac1\alpha}\gamma_n\,|u|)\,
du\right]\right.\\
&&\left.+\frac{\alpha^2e^\alpha}{2(1+\alpha)}\,\E\left[\int_\R \big(\phi'(\widetilde{S}_{n,i}+\frac{n^{-\frac1\alpha}}{\sigma}u)-\phi'(\widetilde{S}_{n,i})\big)
\frac{
u\log(
n^{-\frac1\alpha} \big|u|
)
}
{
|u|^{\alpha +1}\log \gamma_n
}\,
{\bf 1}_{[e,\infty)}(n^{-\frac1\alpha}\gamma_n\,|u|)\,
du\right]\right|\\
&\leq &\frac{\alpha^2e^\alpha}{1+\alpha}\,\|\phi''\|_\infty \frac{n^{-\frac1\alpha}}{\sigma} \int_{0}^{\frac{e{n^\frac1\alpha}}{\gamma_n}}
\frac{du}{u^{\alpha-1}}+\frac{\alpha^2e^\alpha}{1+\alpha}\,\|\phi''\|_\infty \frac{n^{-\frac1\alpha}}{\sigma}\frac{1}{\log\gamma_n}\int_{\frac{en^{\frac{1}{\alpha}}}{\gamma_{n}}}^{n^{\frac{1}{\alpha}}}\frac{|\log(n^{-\frac{1}{\alpha}}u)|}{u^{\alpha-1}}du\\
&&\hskip5cm+\frac{2\alpha^2e^\alpha}{1+\alpha}\,\|\phi'\|_\infty\frac{1}{\log\gamma_n}
\int_{n^{\frac{1}{\alpha}}}^\infty
\frac{|\log (n^{-\frac1\alpha} u) |du}{u^{\alpha}}\\
&=&O(n^{\frac1\alpha-1}\gamma_n^{\alpha-2})+\frac{n^{\frac1\alpha-1}}{\log\gamma_n}\Big(\frac{\alpha^2e^\alpha}{\sigma(1+\alpha)}\,\|\phi''\|_\infty\int_{\frac{e}{\gamma_{n}}}^{1}\frac{|\log(v)|}{v^{\alpha-1}}dv
+\frac{2\alpha^2e^\alpha}{1+\alpha}\,\|\phi'\|_\infty
\int_{1}^\infty
\frac{|\log v |dv}{v^{\alpha}}\Big)\\
&=&O\big(n^{\frac1\alpha-1}(\log\gamma_{n})^{-1}\big).
\end{eqnarray*}
On the other hand, by (\ref{crucialbis}) we have
\begin{eqnarray*}
\left|
\E\big[
(\mathcal{L}^{\alpha,0}\phi)(\widetilde{S}_{n,i})\big]
-\E\big[
(\mathcal{L}^{\alpha,0}\phi)(\widetilde{S}_n)\big]
\right|\leq \frac{4d_{\alpha}\|\phi''\|_\infty}{\alpha(2-\alpha)(\alpha-1)}\E[|X_1|^{2-\alpha}]\sigma^{\alpha-2}\,\gamma_n^{\alpha-2}=O(n^{1-\frac2\alpha}(\log n)^{1-\frac2\alpha}).
\end{eqnarray*}
Putting everything together, we get that
\[
d_W(\widetilde{S}_n,S\alpha S(1))=O\big((\log\gamma_{n})^{-1})=O((\log n)^{-1}\big),
\]
which is the desired conclusion.\qed

\end{document}